\documentclass[reqno]{amsart}     
\usepackage{enumerate,amsmath}     
\usepackage{pstricks,pst-node,pst-coil}   
     
\theoremstyle{plain}      
\newtheorem{thm}{Theorem}     
\newtheorem{theorem}[thm]{Theorem}     
     
\newtheorem{corollary}[thm]{Corollary}

\theoremstyle{remark}      
     
\theoremstyle{definition}

\newtheorem{example}[thm]{Example}     
\newtheorem{examples}[thm]{Examples}     

\def\al{{\alpha}}         
\def\be{{\beta}}

\def\Om{{\Omega}}         
\def\la{{\lambda}}

\def\Si{{\Sigma}}

\def\phi{{\varphi}}

\let\pa\partial     
\let\na\nabla

\DeclareMathAlphabet{\doba}{U}{msb}{m}{n}

\gdef\mN{\doba{N}}

\gdef\mR{\doba{R}}

\gdef\cL{\mathcal{L}}

\def\grad{{\mathop{\rm grad\;}}}     
\def\divv{{\mathop{\rm div}}}

\let\ti\tilde   
\def\eref#1{{\rm (\ref{#1})}}   
\let\<\langle 
\let\>\rangle

\begin{document}     
\title{An obstruction for the mean curvature of a 
conformal immersion $S^n\to \mR^{n+1}$}
\author{Bernd Ammann, Emmanuel Humbert, Mohameden Ould Ahmedou}     
\date{June 2005}     
     
\begin{abstract}
We prove a Pohozaev type identity for non-linear eigenvalue equations
of the Dirac operator on Riemannian spin manifolds with boundary.
As an application, we obtain that the mean curvature $H$ of a conformal
immersion $S^n\to \mR^{n+1}$ satisfies
$\int \pa_X H=0$ where $X$ is a conformal vector field on $S^n$ and where
the integration is carried out with respect to the Euclidean
volume measure of the image.
This identity is analogous to the Kazdan-Warner obstruction that appears in
the problem of prescribing the scalar curvature on $S^n$ 
inside the standard conformal class.
\end{abstract}
   
\maketitle     

Let $(M,g)$ be a compact Riemannian manifold with a conformal 
vector field~$X$. Given a function $s$ on $M$, 
then it is a classical question to ask whether $s$ is the 
scalar curvature of a metric $\ti g$ conformal to $g$.
The determination of 
the set of all such functions $s$ is still open, although several partial 
results are known, in particular there are necessary conditions that
$s$ has to satisfy in order to be a scalar curvature.

On the one hand there are topological obstructions. 
If for example $M$ is spin and has non-vanishing $\hat A$ genus, then the
scalar curvature of any metric on $M$ has either to be negative somewhere
or the Ricci curvature vanishes everywhere on $M$.

However, if one fixes the conformal class $[g]$ as described 
above, there are further obstructions that arise from conformal vector fields.
For example if $M$ is $S^n$ with the standard conformal structure, 
Kazdan and Warner~\cite{kazdan.warner:75} derived a famous obstruction. 
A slightly stronger version of 
this obstruction due to Bourguignon and Ezin 
\cite{bourguignon.ezin:87} is described in the following theorem. 
\begin{theorem}
Let $X$ be a conformal vector field on the compact manifold $(M,g)$.
If $s$ is the scalar curvature of a metric $\ti g=u^{4/(n-2)}g$, then 
  $$\int_M \pa_Xs\;dv_{\ti g}=0 $$
where $dv_{\ti g}=u^{2n\over n-2}\,dv_g$ 
is the volume measure associated to $\ti g$.
\end{theorem}

Tightly related to the Kazdan-Warner obstruction is the Pohozaev identity. 
Let $\Om$ be a star-shaped open set of  $\mR^n$ ($n \in \mN$) with 
smooth boundary. We denote by $\Delta = - \sum_{i=1}^n \pa_{ii}$ 
the Laplacian on $\mR^n$. 
Let $u \in C^2(\bar{\Om})$ be a positive solution of 
$\Delta u = u^{p-1}$ on $\Om$ with
$u_{|\partial \Om} \equiv 0$. The vector field $X=\sum_{i=1}^n x^i\pa_{i}$ is 
conformal. If one uses similar methods
as in the proof of the Kazdan-Warner obstruction,
then one obtains the Pohozaev identity (\cite{pohozaev:65}) which
asserts that:
\begin{equation}\label{eq.pohoz}
    \left(1 -\frac{n}{2}+\frac{n}{p}\right) \int_{\Om} u^p = 
    \frac{1}{2} \int_{\pa \Om} \<\nu,X\> ( \pa_{\nu} u)^2 
\end{equation}
where $\nu$ resp.\ $\pa_{\nu}$ is the outer normal vector
resp.\ the outer normal derivative on~$\pa \Om$.
One among many important consequences of this inequality is that 
no non-trivial solutions exist if $p\geq \frac{2n}{n-2}$. Another 
application is an alternative
proof of the Kazdan-Warner obstructions 
in the case that $(M,g)$ is the sphere with the standard 
conformal structure \cite{druet.robert:99}.

In the present short article, we establish  
a similar identity for the classical 
Dirac operator $D$. We derive this identity on manifolds with boundary
in order to admit future Pohozaev type applications. Then, we will
specialize to compact manifolds without boundary, where 
we will derive a Kazdan-Warner type obstruction
for the mean curvature of a conformal immersion $S^2\to \mR^3$.

Our main theorem is:
\begin{theorem} 
Let $(M,g,\chi)$ be a compact Riemannian spin manifold of dimension 
$n$ with boundary $\partial M$ (possibly equal to 
$\emptyset$) and with Dirac operator $D$.
We assume that there exists a smooth spinor field $\psi$
which satisfies for some $p>1$,
\begin{equation}\label{eq.sol}
  D\psi =H |\psi|^{p-2}\psi, \qquad H\in C^\infty(M).
\end{equation}
Furthermore, we assume that $X$ is a conformal vector field on $M$.
Then, we have the following Pohozaev type identity:
 $$\int_{\pa M} \<\nu\cdot \cL_X \psi,\psi\> 
     = {p-2\over p}\, \int_{\pa M} H |\psi|^p g(X,\nu) 
     +  \left(1-{p-2\over p}\,n\right) \int_M H \beta\, |\psi|^p
     +\frac2p \int_M (\pa_X H) \,|\psi|^p$$
where $\nu$  denotes the outward pointing normal vector along $\pa M$,
and where $\<\,,\,\>$ denotes the real scalar product on spinors.
\end{theorem} 

{\bf Proof:} The flow associated to the conformal vector
field $X$ will be denoted as $\al^t$. If $p$ is in the interior
of $M$, then $\al^t(p)$ exists for times $t$ close to $0$.
For any $t\in \mR$ let $f^t$ be the conformal scaling function of $\al^t$, 
i.e.\ $(d\al^t)_p$ is $f^t(p)$ times an isometry from 
$T_pM$ to $T_{\al^t(p)} M$.
Let $\al^t_*:\Si_pM \to \Si_{\al^t(p)}M$ 
be the spinor identification map as constructed in 
\cite{hitchin:74,hijazi:86,bourguignon.gauduchon:92}. 
In particular, this map has the pointwise properties that 
  $$|\al^t_*(\psi)|=|\psi|$$
and the following transformation formula for 
conformal changes of the metric. 
Let $\phi\in \Gamma(\Si M)$ be a spinor field.
For $t$ close to $0$, we 
then define the map $\al^t_{\#}:\Gamma(\Si M)\to \Gamma(\Si \ti M)$,
$\al^t_{\#}(\phi):=\al^t_*\circ \phi\circ \al^{-t}$, where 
$\ti M$ is $M$ without an open neighborhood of the boundary.

Then  
$$ D\al^t_{\#}\left((f^t)^{-{n-1\over 2}}\psi\right) = \al^t_{\#}((f^t)^{-{n+1\over 2}}D\psi).$$
Now we assume that $\psi$ satisfies \eref{eq.sol}, and we obtain
  $$ D\al^t_{\#}\left((f^t)^{-{n-1\over 2}}\psi\right)
     =\al^t_{\#}\left((f^t)^{-{n+1\over 2}}H |\psi|^{p-2}\psi\right).$$
Deriving with respect to $t$ at $t=0$ yields
\begin{eqnarray}\label{eq.celder}
   -{n-1\over 2} D\beta \psi+  D {d\over dt}|_{t=0}\al^t_\#\psi &=& 
   -{n+1\over 2} H \beta |\psi|^{p-2} \psi
   + H |\psi|^{p-2} {d\over dt}|_{t=0}\al^t_\#\psi \\
   &&+\; (p-2)H \<  {d\over dt}|_{t=0}\al^t_\#\psi,\psi\> |\psi|^{p-4} \psi
      -(\pa_X H)|\psi|^{p-2}\psi. 
\end{eqnarray}
where $\be:={d\over dt}|_{t=0}f^t$.
We reformulate using definition of the \emph{Lie derivative of spinor fields in the direction $X$} 
\cite{bourguignon.gauduchon:92}, i.e.\
\begin{equation}\label{eq.lieder}
  \cL_X(\psi)= - {d\over dt}|_{t=0}\al^t_{\#}(\psi).
\end{equation}
Together with 
$D\beta\psi=\beta D\psi+\na \beta\cdot\psi$ and \eref{eq.sol} one then 
concludes that 
\begin{eqnarray}
{n-1\over 2} \,\na \be \cdot \psi+D\cL_X\psi &=&  H \beta |\psi|^{p-2}\psi 
{}+\, H |\psi|^{p-2} \cL_X\psi\\
 && +\;  (p-2)H \<\cL_X\psi,\psi\>|\psi|^{p-4}\psi + (\pa_X H) |\psi|^{p-2}\psi.
\end{eqnarray}

After multiplication with $\psi$, 
the $\na \be \cdot \psi$-term vanishes, and we obtain
  $$\<D\cL_X \psi,\psi\> =  (p-1) H |\psi|^{p-2} \< \cL_X \psi,\psi\>
    +  H \beta |\psi|^p + (\pa_X H) |\psi|^p.$$
The product rule for the Lie derivative tells us that 
$$|\psi|^{p-2}\<\cL_X\psi,\psi\>= {1\over 2}\,|\psi|^{p-2} \pa_X |\psi|^2 = |\psi|^{p-1} \pa_X |\psi|  = {1\over p}\,\pa_X |\psi|^p.$$
Hence, we obtain
  $$\<D\cL_X \psi,\psi\> = {p-1\over p} \,H\,\pa_X |\psi|^p  
    + H\beta |\psi|^p+ (\pa_X H) |\psi|^p.$$
Strictly speaking, this equation is valid in the interior, but it extends to
the boundary by continuity.
Now, we integrate over $M$.
With partial integration for the Dirac operator one obtains
$$\int_M \<D\cL_X \psi, \psi\> 
    =\int_M \< \cL_X \psi, D\psi\>+ \int_{\pa M} 
    \<\nu\cdot \cL_X \psi,\psi\>
    = \int_M H\underbrace{\< \cL_X \psi,|\psi|^{p-2}\psi\>}_{={1\over p} \pa_X |\psi|^p}+ \int_{\pa M} \<\nu\cdot \cL_X \psi,\psi\>.$$ 
This yields
  $$\int_{\pa M} \<\nu\cdot \cL_X \psi,\psi\> = {p-2\over p} \int_M H \pa_X |\psi|^p  
    + \int_M H \beta |\psi|^p+  \int_M (\pa_X H) |\psi|^p.$$
Using $\divv (H |\psi|^p X) = (\pa_X H) |\psi|^p+ H \pa_X |\psi|^p + H |\psi|^p \divv X$ and $\divv X=n\beta$ we obtain
  $$\int_{\pa M} \<\nu\cdot \cL_X \psi,\psi\> 
     = {p-2\over p}\, \int_{\pa M} H |\psi|^p g(X,\nu) 
     +  \left(1-{p-2\over p}\,n\right) \int_M H \beta |\psi|^p
     +\frac2p \int_M (\pa_X H) |\psi|^p.$$

\begin{examples}\ \\
1.) Let $\Omega$ be domain in $\mR^n$ with smooth boundary, 
let $X=r\pa_r=\sum x^i\pa_i$, and we will assume that $H=\lambda$ is constant.
Then $\beta\equiv 1$ and we obtain
  $$\int_{\pa\Omega} \<\nu\cdot \cL_X \psi,\psi\> = \la \,{p-2\over p} \int_{\pa \Omega} \<X,\nu\>\,|\psi|^p
     +  \la\left(1-{p-2\over p}\,n\right) \int_{\Omega} |\psi|^p.$$
This inequality bears many analogies to equation~\eref{eq.pohoz}. 
In particular, the constant $1-{p-2\over p}$ before the 
integral over $\Omega$ vanishes if $p$ takes the value 
$p=2n/(n-1)$. This value plays the same role in non-linear Dirac 
equations as the value $p=2n/(n-2)$ does for the Laplace operator.\\
2.) If $M$ is a closed manifold and $X$ is a conformal vector field, then for $p=2n/(n-1)$ we obtain
  $$\int_M (\pa_XH)|\psi|^p =0.$$
\end{examples}

\begin{corollary} {\bf [Kazdan-Warner type obstructions]}
Let $f:S^n\to \mR^{n+1}$, $n\geq 2$, 
be a conformal immersion (possibly with branching points
of even order in the case $n=2$).
We denote by $H:S^n \to \mR$ the mean curvature of this
immersion. Then, for any conformal vector field $X$ we have 
 $$\int_{S^n} (\pa_X H) f^*(d\mu) =0$$
where $d\mu$ is the volume element on $f(S^n)$
induced from the euclidean metric on $\mR^3$.
In particular, $\pa_X H$ changes sign. 
\end{corollary} 

The corollary is particularly interesting in dimension $n=2$. If 
$f:S^2\to \mR^3$
is any immersion, then after possibly composing with a diffeomorphism
$S^2\to S^2$, we can assume that $f$ is conformal.

The corollary is analogous to results in
\cite{kazdan.warner:75},  \cite{bourguignon.ezin:87} and 
\cite{druet.robert:99}.

{\bf Proof:}  Let $\psi$ be parallel spinor on $\mR^{n+1}$. Then, as
proven in \cite{kusner.schmitt:p96,baer:98,friedrich:98}, 
the restriction of $\psi$ on $\Si$ satisfies equation
(\ref{eq.sol}) with $p=2n/(n-1)$, and $|\psi|^p\,d\nu=f^*(d\mu)$ 
where $d\nu$ is the standard volume element on~$S^n$. 
Since this equation is conformally invariant
we obtain a solution of (\ref{eq.sol}) on $S^n$ equipped with the standard
metric. The corollary then immediately follows from example (2)
above. 

\begin{example} 
Let $x_3:S^2\to \mR$ be the third component of the
standard inclusion.
One shows that $X:=\grad x_3$ is a conformal vector
field on $S^2$, where the gradient is taken with respect to the 
standard metric on $S^2$. Then for any $C\in \mR$ one has
$\pa_X(x_3+C)=g(\grad x_3,\grad x_3)\geq 0$.
Hence $x_3+C:S^2\to \mR$ is not the mean curvature of a conformal immersion.  
\end{example}

 
\long\def\komment#1{}

             
           
\vspace{1cm}               
Authors' addresses:               
\nopagebreak   
\vspace{5mm}\\   
\parskip0ex   
             
\vtop{   
\hsize=8cm\noindent   
\obeylines               
Bernd Ammann and Emmanuel Humbert,               
Institut \'Elie Cartan BP 239              
Universit\'e de Nancy 1               
54506 Vandoeuvre-l\`es -Nancy Cedex               
France                           
\vspace{0.5cm}               

\vtop{   
\hsize=8cm\noindent   
\obeylines               
M. Ould Ahmedou
Mathematisches Institut
der Universit\"at T\"ubingen
Auf der Morgenstelle 10
72076 T\"ubingen
Germany
\vspace{0.5cm}               
}

E-Mail:               

{\tt ammann at iecn.u-nancy.fr,  humbert at iecn.u-nancy.fr and ahmedou at analysis.mathematik.uni-tuebingen.de}

WWW:               

{\tt http://www.berndammann.de/publications}               

}


\end{document}